\documentclass{amsart}

\newtheorem{theorem}{Theorem}[section]

\theoremstyle{definition}

\newtheorem{proposition}[theorem]{Proposition}

\theoremstyle{remark}
\newtheorem{remark}[theorem]{Remark}

\begin{document}

\title[Weighted polynomial multiple ergodic averages]{Convergence of weighted polynomial multiple ergodic averages}

\author{Qing Chu}

\address{Universit\'{e} Paris-Est, Laboratoire d'Analyse et de math\'{e}matiques appliqu\'{e}es, UMR
                                               ´
CNRS 8050, 5 bd Descartes, 77454 Marne la Vall\'{e}e Cedex 2, France}

\email{qing.chu@univ-mlv.fr}


\date{\today}%


\begin{abstract}
In this article we study weighted polynomial multiple ergodic averages. A
sequence of weights is called universally good if any polynomial
multiple ergodic average with this sequence of weights converges in
$L^{2}$. We find a necessary condition and show that for any bounded
measurable function $\phi$ on an ergodic system, the sequence
$\phi(T^{n}x)$ is universally good for almost every $x$.  The linear
case was covered by Host and Kra.
\end{abstract}

\maketitle

\section{introduction}
In his innovative proof of Szemer\'edi's Theorem via ergodic theory,
Furstenberg introduced certain multiple ergodic averages.  There
have been many results on these and other nonconventional ergodic
averages, including the multiple ergodic theorems of Host and
Kra~\cite{hk1}, \cite{hk2}, Leibman~\cite{Lei}, Ziegler~\cite{Z}.
Recently Host and Kra studied weighted ergodic theorems for multiple
averages along arithmetic progressions. We give a generalization of
this result for polynomial averages, showing:

\begin{theorem}\label{main0}
Let $(Y,\nu,S)$ be an ergodic system and $\phi\in L^\infty(\nu)$.
Then there exists $Y_0\subset Y$ with $\nu(Y_0)=1$ such that, for
every $y_0\in Y_0$, every system $(X,\mu,T)$, every $r\geq 1$, all
integer polynomials $p_1,\dots,p_r$ and all functions
$f_1,\dots,f_r\in L^\infty(\mu)$, the averages
$$
\frac 1N\sum_{n=0}^{N-1}\phi(S^ny_0)
T^{p_{1}(n)}f_{1}\cdots T^{p_{r}(n)}f_{r}
$$
converge in $L^2(\mu)$ .
\end{theorem}
Throughout this article, by {\em integer polynomial} we mean a
polynomial all of whose coefficients are integers.

The case of $p_{i}(n)= in$ was proved by Host and Kra \cite{hk}.
\begin{remark}
One may wonder why in the theorem $S^n$ is not replaced by $S^{p(n)}$ for some integer
polynomial $p(n)$. In fact, the latter would be much harder to
prove, even in the simplest case that $r= 1$ and $p_1(n)= n$. This problem is equivalent to showing the convergence of the averages of $\phi(S^{p(n)} y_0)e^{2\pi int}$ for all $t\in \mathbb{T}$ and all $y_0 \in Y_0$, where $Y_0$ does not depend on $t$. But this problem reduces to a question of almost everywhere convergence of multiple ergodic averages along polynomials, and this question is out of reach for the moment.
\end{remark}

Note that the set $Y_0$ does not depend on $X$ or on $f_{i},
i=1,\dots,r$. We say that for every $y_0\in Y_0$, the sequence
$\phi(S^ny_0)$ is \emph{universally good for the convergence in the mean
of polynomial multiple ergodic averages}.

For $r=1$ and $p(n)=n$, the result follows immediately from the
classical Wiener-Winter ergodic Theorem \cite{ww} and a corollary of
the Spectral Theorem. We follow a similar strategy, generalizing the
proof in~\cite{hk} along arithmetic progressions to polynomial
progressions, but need to address some deeper technical issues.

We first recall some definitions, see ~\cite{BHK} and ~\cite{hk} for
details. Let $G$ be a $k$-step nilpotent Lie group and $\Gamma
\subset G$ be a discrete, cocompact subgroup of $G$.  The compact
manifold $X = G/\Gamma$ is called a \emph{$k$-step nilmanifold}.
The Haar measure $\mu$ of $X$ is the unique probability measure
invariant under the left translations $x \mapsto g x$ of $G$ on $X$.
Letting $T$ denote left multiplication by a fixed element $\alpha
\in G$, we call $(X, \mu, T )$ a \emph{$k$-step nilsystem}.  Let $f
\colon X \to C$ be a continuous function, $x_{0}\in X$, then the
sequence $(f (\alpha^{n} x_0 ) : n \in \mathbb{Z})$ is called a
\emph{basic $k$-step nilsequence}. The family of basic $k$-step nilsequences forms a subalgebra of $l^\infty$. Under the uniform norm $|\!|\cdot|\!|_{\infty}$ of $l^{\infty}$, we call a
uniform limit of basic $k$-step nilsequences a \emph{k-step nilsequence}.

\

The proof of Theorem \ref{main0} is broken down into two pieces. First we give a convergence criterion for weighted polynomial multiple ergodic averages.

\begin{theorem}[Convergence Criterion for Weighted Averages]\label{main}
For any $r,b\in\mathbb{N}$, there exists an integer $K\geq 1$ with
the following property: for any bounded sequence $\mathbf{c}=
(c_{n}: n\in \mathbb{Z})$, if the averages
\begin{equation*}
\frac{1}{N}\sum_{n=0}^{N-1}c_{n}d_{n}
\end{equation*}
converge as $N\rightarrow \infty$ for every $K$-step nilsequence
$\textbf{d}= (d_{n}: n\in\mathbb{Z})$. Then for every system $(X,
\mu,T)$, all $f_{1},\dots,f_{r}\in L^{\infty}(X)$, and all integer
polynomials $p_{1},\dots,p_{r}$ of degree $\leq b$, the averages
\begin{align}\label{main1}
\frac{1}{N}\sum_{n=0}^{N-1}c_{n}T^{p_{1}(n)}f_{1}\cdot
T^{p_{2}(n)}f_{2}\cdots T^{p_{r}(n)}f_{r}
\end{align}
converge in $L^{2}(X)$.
\end{theorem}

The bulk of this paper is devoted to the proof of this theorem. Then
our main result follows from the following Generalized
Wiener-Wintner Theorem proved by Host and Kra in ~\cite{hk}. The
case of a polynomial version of the Wiener-Wintner theorem was proved by
Lesigne (~\cite{lesi1}, ~\cite{lesi}).

\begin{theorem}[Generalized Wiener-Wintner Theorem \cite{hk}]\label{GWW}
Let $(X, \mu, T )$ be an ergodic system and $\phi$ be a bounded
measurable function on $X$. Then there exists $X_0 \subset X$ with
$\mu (X_0) = 1$ such that for every $x \in X_0$, the averages
$$
\frac{1}{N}\sum_{n=0}^{N-1}\phi(T^{n}x)b_{n}
$$
converge as $N\rightarrow\infty$ for every $x\in X_0$ and every
nilsequence $b = (b_{n} : n \in \mathbb{Z})$.
\end{theorem}

While nilsequences do not appear in the statement of
Theorem~\ref{main0}, they are used as tools in its proof.  Both
Theorems~\ref{main} and~\ref{GWW} are of interest on their own, as
results on nilsequences.

\subsection*{Acknowledgement}
We thank the referee for very useful remarks and the simplification of the proof. 

\section{proof of theorem \ref{main}}
Using a standard ergodic decomposition argument, it suffices to prove Theorem \ref{main} for ergodic systems.

\subsection{}
We first remind the reader of a definition from Leibman's paper \cite{Lei2}. We call a sequence
$\{g(n)\}_{n\in \mathbb{Z}}$ with values in a nilpotent group $G$ a
\emph{polynomial sequence}, if $g(n)$ is of the form $g(n)=
a_1^{p_1(n)}\ldots a_m^{p_m(n)}$, where $a_1,\dots,a_m \in G$ and
$p_1,\dots,p_m$ are polynomials taking integer values on the
integers.

Before stating the next proposition, we explain briefly its meaning: we can view the sequence of values of a continuous function
along a polynomial sequence on a nilmanifold as the sequence of values of some
other continuous function along an ordinary orbit of some other nilsystem.

\begin{proposition}\label{nil}
Let $(X= G/\Gamma, T)$ be a nilsystem, $x_{0}\in X$, p be an
\emph{integer polynomial}, and $f\in \mathcal{C}(X)$. Then there
exists a nilsystem $(Y , S)$, $y_{0}\in Y$, $h\in \mathcal{C}(Y)$,
such that $f(T^{p(n)}x_{0})= h(S^{n}y_{0})$ for every $n$.
\end{proposition}

\begin{proof}
Let $(X= G/\Gamma, T)$ be a nilsystem. Suppose $Tx:= \alpha x$, for
some $\alpha\in G$. Then $T^{p(n)}x= \alpha ^{p(n)}x$. Let $g(n):=
\alpha ^{p(n)}$, then $g$ is a polynomial sequence in $G$. Let $\pi:
G\rightarrow X$ be the factorization mapping. We will assume that
$x_{0}= \pi (\textbf{1}_{G})$; otherwise if $x_0= \pi(\gamma),
\gamma \in G$, we write $g(n)x_0= g(n)\gamma \pi(\textbf{1}_{G})$,
and replace $g(n)$ by $g(n)\gamma$.

Now we have a nilpotent Lie group $G$, a discrete cocompact subgroup
$\Gamma$ and a polynomial sequence $g$ in $G$. By Proposition 3.14
in Leibman's paper \cite{Lei2}, there exist a nilpotent Lie group
$\widetilde{G}$, a discrete cocompact subgroup $\widetilde{\Gamma}$,
an epimorphism $\eta: \widetilde{G}\rightarrow G$ with
$\eta(\widetilde{\Gamma})\subseteq \Gamma$, a unipotent automorphism
$\widetilde{\tau}$ of $\widetilde{G}$ with
$\widetilde{\tau}(\widetilde{\Gamma})= \widetilde{\Gamma}$, and an
element $\tilde{c}\in \widetilde{G}$ such that $$g(n)=
\eta(\widetilde{\tau}^{n}(\tilde{c})), n\in \mathbb{Z}.$$

Let $\widetilde{X}= \widetilde{G}/\widetilde{\Gamma}$ and let
$\widetilde{\pi}: \widetilde{G}\rightarrow \widetilde{X}$ be the
factorization mapping.

The epimorphism $\eta: \widetilde{G}\rightarrow G$ factors to a map
$\widetilde{X}\rightarrow X$, we also denote it by $\eta$, which is
onto and satisfies,
$$
\pi\circ\eta= \eta\circ\widetilde{\pi} .
$$

The map $\widetilde{\tau}$ induces a homomorphism $\widetilde X\to\widetilde X $, which we
also denote it by $\widetilde{\tau}$. It satisfies,
$$
\widetilde{\tau}\circ\widetilde\pi=\widetilde\pi\circ\widetilde{\tau} .
$$

Let $\widetilde {x}_{0}=
\widetilde{\pi}(\textbf{1}_{\widetilde{G}})$, then
$$
\eta(\widetilde{\tau}^{n}(\tilde{c}\ \widetilde {x}_{0}))=
g(n)x_{0},\ n\in \mathbb{Z} .
$$

Let $\widehat{G}$ be the extension of $\widetilde{G}$ by
$\widetilde{\tau}$, then $\widehat{G}$ is a nilpotent Lie group (see
Proposition 3.9 in \cite{Lei2}). Let $\widehat{\tau}$ be the element
in $\widehat{G}$ representing $\widetilde{\tau}$, so that
$\widetilde{\tau}(\widetilde{\alpha})=
\widehat{\tau}\widetilde{\alpha}\widehat{\tau}^{-1}$ for any
$\widetilde{\alpha}\in \widetilde{G}$, the multiplication of
$\widehat{G}$ is given by this formula. We have $\widehat{G}= \{\widetilde{g}\widehat{\tau}^{n}\colon \widetilde{g}\in \widetilde{G}, n\in \mathbb{Z}\}$ and $\widetilde{G}$ is open in $\widehat{G}$.

Let $\widehat{\Gamma}$ be the subgroup of $\widehat G$ spanned by
$\widetilde\Gamma$ and $\widehat\tau$. As $\widetilde{\tau}(\widetilde{\Gamma})= \widetilde{\Gamma}$, we have $\widehat{\Gamma}= \{\widetilde{\gamma}{\widehat{\tau}}^{n}\colon \widetilde{\gamma}\in \widetilde{\Gamma}, n\in\mathbb{Z}\}$ and 
$\widehat{\Gamma}\cap \widetilde{G}= \widetilde{\Gamma}$. By the definition of
the relative topology, $\widehat{\Gamma}$ is a discrete
subgroup of $\widehat{G}$. Moreover, 
$\widehat{G} / \widehat{\Gamma}= (\widetilde{G}\widehat{\Gamma}/\widehat{\Gamma})$ can be identified with  $\widetilde{G}/(\widehat{\Gamma}\cap\widetilde{G})= \widetilde{G}/\widetilde{\Gamma}= \widetilde X $. We write $\widehat\pi\colon\widehat G\to\widetilde X$ for the quotient map.

Let $\widetilde x\in\widetilde X$ and $\tilde g\in\widetilde G$
with $\widetilde \pi(\tilde g)=\widetilde x$. We have
$$
\widetilde \tau(\widetilde x)=\widetilde \pi(\widetilde
\tau(\tilde g)) =\widehat\pi(\widehat\tau \tilde g
\widetilde\tau^{-1})=\widehat\pi(\widehat \tau \tilde g)
=\widehat\tau\widetilde x ,
$$
 because $\widehat\tau^{-1}\in\widehat\Gamma$. So for every $n$,
$$
g(n)x_{0}=
\eta(\widetilde{\tau}^{n}(\tilde{c}\ \widetilde {x}_{0}))=
\eta(\widehat{\tau^{n}}(\tilde{c}\ \widetilde {x}_{0})) .
$$

Let $Y=(\widehat{G}/\widehat{\Gamma}, S)= (\widetilde{X}, S)$,
$S\widetilde{x}= \widehat{\tau}\widetilde{x}$, and let $h= f\circ
\eta$, and $y_{0}= \tilde{c}\ \widetilde {x}_{0}$. This system and
this function satisfy the announced properties.
\end{proof}

\subsection{}\label{nota}
We recall a few properties of the seminorms and the factors
introduced in \cite{hk1}. Let $(X, \mu, T )$ be an ergodic system.
For an integer $k\geq 0$, we write $X^{[k]}= X^{2^{k}}$ and $T^{[k]}
\colon X^{[k]}\rightarrow X^{[k]}$ for the map $T\times T\times
\ldots \times T$, taken $2^{k}$ times. We define by induction a probability
measure $\mu^{[k]}$ on $X^{[k]}$ that is invariant under $T^{[k]}$. Set $\mu^{[0]}=\mu$. For $k\geq 0$, let
$\mathcal{I}^{[k]}$ be the $\sigma$-algebra of $T^{[k]}$-invariant
subsets of $X^{[k]}$. Then $\mu^{[k+1]}$ is the relatively
independent product of $\mu^{[k]}$ over $\mathcal{I}^{k}$, which
means for $F, F' \in L^{\infty}(X^{[k]})$,
$$\int_{X^{[k+1]}}F\otimes F'\ d\mu^{[k+1]}= \int_{X^{[k]}}
\mathbb{E}(F|\mathcal{I}^{[k]})\mathbb{E}(F'|\mathcal{I}^{[k]}) \
d\mu^{[k]}.$$

For a bounded measurable function $f$, we define

\begin{align}\label{semi}
|\!|\!|f|\!|\!|_{k} = \left(\int_{X^{[k]}}\prod_{\varepsilon \in \{0,
1\}^{k}} C^{|\varepsilon|}f(x_{\varepsilon})\
d\mu^{[k]}(x)\right)^{1/2^{k}},
\end{align}
where $C\colon \mathbb C\rightarrow \mathbb C$ is the conjugacy map
$z\mapsto \overline{z}$, $\varepsilon =
\varepsilon_1\varepsilon_2\dots\varepsilon_k$ with $\varepsilon_i
\in \{0,1\}$ and
$|\varepsilon|=\varepsilon_1+\varepsilon_2+\dots+\varepsilon_k.$ It
is shown in \cite{hk1} that for every $k \geq 1$, $|\!|\!| \cdot
|\!|\!|_{k}$ is a seminorm on $L^{\infty}(\mu)$.

Moreover, for every $k \geq 2$, $X$ admits a factor $Z_{k-1}$ such
that, for every $f \in L^{\infty} (\mu)$, $|\!|\!|f|\!|\!|_{k} = 0 $
if and only if $\mathbb{E}(f|Z_{k-1}) = 0$. One of the main results
of \cite{hk1} is that, for every $ k$, $Z_{k}$ is an inverse limit
of $k$-step nilsystems. We call this result the Structure
Theorem.

Let $(Z=Z_1(X), m, T)$ be the Kronecker factor of $(X, \mu, T)$. For
$s\in Z$, we define a measure $\mu_s$ on $X\times X$ by
$$\int_{X\times X}f(x)f'(x')\ d\mu_{s}(x,x')= \int_{Z} \mathbb{E}(f|Z)(z)
\cdot \mathbb{E}(f'|Z)(sz)\ d m(z).$$ 

For every $s\in Z$ the measure
$\mu_s$ is invariant under $T\times T$ and is ergodic for $m$-almost
every $s$. The ergodic decomposition of $\mu\times\mu$ under
$T\times T$ is $$\mu\times\mu = \int_{Z}\mu_{s}\ d m(s).$$

For each $s\in Z$ such that $(X\times X, \mu_s, T\times T)$ is
ergodic, and for each integer $k\geq 1$, a measure $(\mu_s)^{[k]}$ on
$(X\times X)^{[k]}$ can be defined in the same way as $\mu^{[k]}$. Furthermore, a seminorm
$|\!|\!|\cdot|\!|\!|_{k,s}$ on $L^{\infty}(\mu_s)$ can be
associated to this measure in the same way as the seminorm
$|\!|\!|\cdot|\!|\!|_k$ is associated to $\mu^{[k]}.$ It follows
from the definition (\ref{semi}) that for every $f\in
L^{\infty}(\mu)$,

\begin{align}\label{relat}
|\!|\!|f|\!|\!|_{k+1}^{2^{k+1}}=\int_{Z} |\!|\!|f\otimes
\bar{f}|\!|\!|_{k,s}^{2^{k}}\ d m(s).
\end{align}

\subsection{}
We return to the proof of theorem \ref{main}. We may assume that the
polynomials $p_{1},\dots,p_{r}$ are nonconstant and essentially
distinct, that is $p_{i}-p_{j}\neq$ constant for $i\neq j$.

The following theorem will be proved in the next section.
\begin{theorem}\label{chara}
For any $r,b\in\mathbb{N}$, there exists $k\in \mathbb{N}$, such that
for any nonconstant essentially distinct polynomials
$p_{1},\dots,p_{r}: \mathbb{Z}\rightarrow \mathbb{Z}$ of degree
$\leq b$, for every ergodic system $(X, \mu,T)$, every
$f_{1},\dots,f_{r}\in L^{\infty}(X)$ with
$|\!|\!|f_{1}|\!|\!|_{k}=0$, and any bounded sequence $\mathbf{c}=
(c_{n}: n\in \mathbb{Z})$, one has
\begin{align}\label{chara1}
\lim_{N\rightarrow\infty}\left\|\frac{1}{N}\sum_{n=0}^{N-1}c_{n}
T^{p_{1}(n)}f_{1}\cdots T^{p_{r}(n)}f_{r}\right\|_{L^{2}(X)}=0.
\end{align}
\end{theorem}

\subsection{}
Now we give the proof of Theorem \ref{main}.
\begin{proof}[\textbf{Proof of Theorem \ref{main} from Proposition \ref{nil} and Theorem \ref{chara}}]
The proof is exactly the same as the proof of Theorem 2.24 in
\cite{hk}. For any $r,b \in \mathbb{N}$, let $k\in \mathbb{N}$ be
the integer in Theorem \ref{chara}, let $Z_{k-1}$ be the
$(k-1)$-th factor of $(X,\mu,T)$, as given by the
Structure Theorem. By definition, if $\mathbb{E}(f_{1}|Z_{k-1})= 0$,
then $|\!|\!|f_{1}|\!|\!|_{k}= 0$, and by Theorem \ref{chara}, the
averages (\ref{chara1}) converge to zero in $L^{2}(X)$. We say that
the factor $Z_{k-1}$ is the characteristic for the convergence of
these averages. Therefore, it suffices to prove the result when the
functions are measurable with respect to the factor $Z_{k-1}$.

Since $Z_{k-1}$ is an inverse limit of $(k-1)$-step
nilsystems by density, we can assume that $(X,\mu,T)$ is a
$(k-1)$-step nilsystem and that the functions
$f_{1},\dots,f_{r}$ are continuous.

But in this case, by Proposition \ref{nil}, for every $x\in X$, and
all polynomials $p_{1},\dots,p_{r}$, there exist nilsystems
$(Y_{1},S_{1}),\dots,(Y_{r},S_{r})$, $y_{i}\in Y_{i}$, and $g_{i}\in
\mathcal{C}(Y_{i})$, such that $f_{i}(T^{p_{i}(n)}x)=
g_{i}(S_{i}^{n}y_{i})$, $i=1,\dots,r$.

Let $K$ be the maximal order of the nilsystems $(Y_{i},
S_{i}),i=1,\dots,r$. Then the system $(Y= Y_{1}\times\cdots\times
Y_{r}, S= S_{1}\times\cdots\times S_{r})$ is a $K$-step nilsystem. Let
$g: Y_{1}\times\cdots\times Y_{r}\rightarrow \mathbb{R}$ be given by
$g(y)= g (y_{1},\dots,y_{r})= g_{1}(y_{1})\cdot\ldots\cdot
g_{r}(y_{r})$. So the sequence
\begin{align*}
&\ \ \{f_{1}(T^{p_{1}(n)}x)\cdot f_{2}(T^{p_{2}(n)}x)\cdot\ldots
\cdot f_{r}(T^{p_{r}(n)}x)\}_{n\in \mathbb{Z}}\\&=
\{g_{1}(S_{1}^{n}y_{1})\cdot g_{2}(S_{2}^{n}y_{2})\cdot\ldots\cdot
g_{r}(S_{r}^{n}y_{r})\}_{n\in \mathbb{Z}}\\&= \{g(S^{n}y)\}_{n\in
\mathbb{Z}}
\end{align*}
is a $K$-step nilsequence and by hypothesis, the averages
(\ref{main1}) converge for every $x\in X$.
\end{proof}

\section{proof of theorem \ref{chara}}
Note that our goal is very similar to the main result in
Leibman's paper \cite{Lei}, the only difference being that in our case
we are dealing with the weighted averages. In fact, we can deduce Theorem \ref{chara}
by following very closely the arguments of Leibman in his paper to
cover our weighted case with some modifications. But here we adopt
another way that allows us to deduce it directly from
Leibman's result.

\begin{proof}[\textbf{Proof of Theorem \ref{chara}}]
We prove the result by several steps.
\

(i) The following assertion follows immediately from Theorem 3 in Leibman's paper
\cite{Lei}:

Let ${r,b} \in \mathbb{N}$ be fixed. There exists an integer $k= k(r, b)$ such that for every family of nonconstant essentially distinct polynomials $p_1,\dots,p_{r} : \mathbb{Z}^{2}\rightarrow \mathbb{Z}$ of degree $\leq b$, we have: for every ergodic system $(X, \mathcal{B}, \mu, T)$, and every $f_1,\dots, f_r\in L^{\infty}(X)$ with $|\!|\!|f_1|\!|\!|_{k} = 0$, one has
$$\lim_{N\rightarrow \infty}\left \|\frac{1}{N^2}\sum_{0\leq {m,n} <N }T^{p_1 (m,n)}f_1 \cdots T^{p_{r}(m,n)}f_{r}\right\|_{L^{2}(X)} = 0.$$

(ii) We use the notation of Section \ref{nota}. It follows from (\ref{relat}) that $|\!|\!|f_1|\!|\!|_{k+1} =
0$ implies $|\!|\!|f_1\otimes \bar{f_1}|\!|\!|_{k, s}= 0$ for
$m$-almost every $s$. Using the previous result to the ergodic system $(X\times
X, \mu_s, T\times T)$, one gets: if $|\!|\!|f_1|\!|\!|_{k+1} = 0$,
then
\begin{align*}
&\ \ \lim_{N\rightarrow \infty}\frac{1}{N^{2}}\sum_{0\leq m, n <N}\left|\int
T^{p_1(m,n)}f_1\cdots T^{p_r(m,n)}f_r\ d\mu \right|^2
\\&=\lim_{N\rightarrow \infty}\frac{1}{N^{2}}\sum_{0\leq m, n
<N}\left|\int_{Z}\int_{X\times X} (T\times T)^{p_1(m,n)}f_1\otimes
\bar{f_1}\cdots (T\times T)^{p_r(m,n)}f_r\otimes\bar{f_r}\
d\mu_{s}dm(s)\right|\\ &\leq \int_Z\lim_{N\rightarrow
\infty}\int_{X\times X}\left|\frac{1}{N^{2}}\sum_{0\leq m, n <N} (T\times
T)^{p_1(m,n)}f_1\otimes \bar{f_1}\cdots (T\times
T)^{p_r(m,n)}f_r\otimes\bar{f_r}\right|d\mu_{s}dm(s)\\ &\leq \int_Z \lim_{N\rightarrow
\infty}\left\|\frac{1}{N^{2}}\sum_{0\leq m, n <N} (T\times
T)^{p_1(m,n)}f_1\otimes \bar{f_1}\cdots (T\times
T)^{p_r(m,n)}f_r\otimes\bar{f_r}\right\|_{L^2(\mu_s)}d m(s)= 0.
\end{align*}

(iii) Expanding the square, one sees that
$$\left\|\frac{1}{N}\sum_{n=0}^{N-1}c_nT^{p_1(n)}f_1\cdots T^{p_r(n)}f_r\right\|_{L^{2}}^2$$
is equal to
$$\frac{1}{N^2}\sum_{0\leq m,n <N}c_n\bar{c}_{m}\int T^{p_1(n)}f_1\cdots T^{p_r(n)}f_r\ 
T^{p_1(m)}\bar{f_1}\cdots T^{p_r(m)}\bar{f_r}\ d\mu,$$
which is less than a constant times

\begin{align}\label{last}
\frac{1}{N^2}\sum_{0\leq m,n <N}\left|\int T^{p_1(n)}f_1\cdots
T^{p_r(n)}f_r\ T^{p_1(m)}\bar{f_1}\cdots
T^{p_r(m)}\bar{f_r}\ d\mu\right|,
\end{align}
since by our assumption the sequence $c_n$ is bounded. Notice that
$p_1(n),\dots,p_r(n),p_1(m),\dots,p_r(m)$ is a family of $2r$
essentially distinct polynomials of degree at most $b$.

(iv) Combining the previous parts we see that if
$|\!|\!|f_1|\!|\!|_{k(2r,b)+1}= 0$, where $k(r,b)$ is
defined in part (i), then for every family of nonconstant 
essentially distinct polynomials $p_{1},\dots,p_{r}:
\mathbb{Z}\rightarrow \mathbb{Z}$ of degree $\leq b$, we have that the average
(\ref{last}) converges to $0$. This proves Theorem \ref{chara}.

\end{proof}

\bibliographystyle{amsplain}

\begin{thebibliography}{99}
\bibitem{BHK}
V.~Bergelson, B.~Host and B.~Kra, with an Appendix by I.~Ruzsa.
\textit{Multiple recurrence and nilsequences.}  Inventiones Math. \textbf{
160} (2005), 261-303.
\bibitem{hk1} B. Host and B. Kra. \textit{Nonconventional ergodic averages and nilmanifolds.} Ann. Math. \textbf{161}
      (2005), 397-488.
\bibitem{hk2} B. Host and B. Kra. \textit{Convergence of polynomial ergodic averages.} Isr. J. Math. \textbf{149} (2005), 1-19.
\bibitem{hk} B. Host and B. Kra. \textit{Uniformity seminorms on $l^{\infty}$ and applications.} Available at
http://arxiv.org/abs/0711.3637v1.
\bibitem{Lei} A. Leibman. \textit{Convergence of multiple ergodic averages along polynomials of several
variables.} Isr. J. Math. \textbf{146} (2005), 303-315.
\bibitem{Lei2} A. Leibman. \textit{Pointwise convergence of ergodic averages for polynomial sequence of translation on a
nilmanifold.} Erg. Th. and Dyn. Sys. \textbf{25} (2005), 201-213.
\bibitem{lesi1} E. Lesigne. \textit{Un th\'{e}or\`{e}me de disjonction de syst\`{e}mes dynamiques et une g\'{e}n\'{e}ralisation du th\'{e}or\`{e}me ergodique de Wiener-Wintner.} Erg. Th. and Dyn. Sys. \textbf{10} (1990), 513-521.
\bibitem{lesi} E. Lesigne. \textit{Spectre quasi-discret et th\'{e}or\`{e}me ergodique de Wiener-Wintner pour les
        polyn\^{o}mes.} Erg. Th. and Dyn. Sys. \textbf{13} (1993), 767-784.
\bibitem{ww} N. Wiener and A. Wintner. \textit{Harmonic analysis and ergodic theory.} Amer. J. Math. \textbf{63}
     (1941), 415-426.
\bibitem{Z} T. Ziegler. \textit{A non-conventional ergodic theorem for a nilsystem.} Erg. Th. and Dyn. Sys. \textbf{25} (2005), 1357-1370.
\end{thebibliography}

\end{document}